# Precise Definition and Analyticity of $\zeta^{(r)}(s,\alpha)$

Vivek V. Rane

C/o Institute of Science,

Madame Cama Road,

Mumbai-400 032, India .

e-mail v_v_rane @yahoo.co.in

Abstract : S. Ramanujan[1] was aware of the power series expression in α of $\zeta(s,1-\alpha)$ for complex s and for $0 \leq \alpha < 1$, which he did not explore very far . Author[2] had derived the power series expression in α of $\zeta(s,1+\alpha)$ for complex α with $|\alpha|<1$ and had shown in author[3] that the power series of $\zeta(-n,\alpha)$ for integral $n \geq 0$, is a polynomial in α .On this backdrop , we give here the precise definition of $\zeta^{(r)}(s,\alpha) = \frac{\partial^r}{\partial s^r}\zeta(s,\alpha)$ in terms of power series in α for complex s and α . We also discuss the analyticity of $\zeta^{(r)}(s,\alpha)$ as functions of s and α and show that $\frac{\partial}{\partial \alpha}\frac{\partial^r}{\partial s^r}\zeta(s,\alpha) = \frac{\partial^r}{\partial s^r}\frac{\partial}{\partial \alpha}\zeta(s,\alpha)$. We compute $\frac{\partial}{\partial \alpha}\frac{\partial^r}{\partial s^r}\zeta(s,\alpha)$ for $s=0$ and $s=1$. We discuss the power series in s of $\zeta(s+1,\alpha)$ for complex α .

# Precise Definition and Analyticity of $\zeta^{(r)}(s,\alpha)$

## Vivek V. Rane

For $0 < \alpha \leq 1$ and for a complex variable $s = \sigma + it$, where $\sigma$ and t are real, let $\zeta(s,\alpha)$ be the Hurwitz's zeta function defined by $\zeta(s,\alpha) = \sum_{n \geq 0} (n+\alpha)^{-s}$ for Re s>1 ; and its analytic continuation .

S. Ramanujan[1] was aware of the power series expression for $\zeta(s, 1-\alpha)$ for $0 \leq \alpha < 1$. However , he did not explore it very far . In author [2] , the power series expression in α for $\zeta(s, 1+\alpha)$ has been derived for complex α with $|\alpha| < 1$; and in author [3] , the power series in α for $\zeta(-n, \alpha)$ for integral $n \geq 0$, has been shown to be a polynomial . The theory of $\zeta^{(r)}(s,\alpha) = \frac{\partial^r}{\partial s^r} \zeta(s,\alpha)$ for integral $r \geq 0$ , as a function of α , can be seen to be a missing link among Riemann zeta function , Euler's gamma function , Bernoulli polynomials, Euler's polynomials, Gauss' ψ- function , Dirichlet L-series , Euler-Zagier's zeta function , Tornheim zeta function etcetera. On this backdrop , we show the following .

Let $k \geq 1$ be an arbitrarily fixed integer . In what follows , empty sum denotes 0 and empty product denotes 1 . For integral $r \geq 0$ , $\zeta^{(r)}(s,\alpha)$ shall stand for $\frac{\partial^r}{\partial s^r} \zeta(s,\alpha)$ and in general , the superscript (r) shall denote the r-th order partial derivative with respect to s .

Theorem : 1) We have for $|\alpha| < k$ ,

$$\zeta(s,\alpha) = \sum_{0 \leq n \leq k-1} (n+\alpha)^{-s} + \sum_{n \geq 0} c_{n,k}(s) \frac{(-\alpha)^n}{n!} ,$$

where
$$c_{n,k}(s) = s(s+1)\ldots\ldots\ldots(s+n-1)\zeta_k(s+n)$$

with
$$\zeta_k(s) = \zeta(s) - \sum_{n \leq k-1} n^{-s}.$$

2) $\zeta(s,\alpha) - \zeta(s)$ is an entire function of s and we have for $r \geq 0$,

$$\zeta^{(r)}(s,\alpha) - \zeta_k^{(r)}(s) = \sum_{0 \leq n \leq k-1} \frac{\partial^r}{\partial s^r}(n+\alpha)^{-s} + \sum_{n \geq 1} \frac{\partial^r}{\partial s^r}(s(s+1)\ldots\ldots(s+n-1)\zeta_k(s+n))\frac{(-\alpha)^n}{n!}$$

for $|\alpha| < k$ except possibly for $\alpha = 0,-1,-2,\ldots\ldots\ldots$.

3) For $r \geq 0$, $\zeta^{(r)}(s,\alpha) - \zeta^{(r)}(s)$ is an analytic function of α except possibly for

$\alpha = 0,-1,-2,\ldots\ldots$, and we have for $\alpha \neq 0,-1,-2,\ldots\ldots$,

$$\frac{\partial}{\partial \alpha}\frac{\partial^r}{\partial s^r}\zeta(s,\alpha) = \frac{\partial^r}{\partial s^r}\frac{\partial}{\partial \alpha}\zeta(s,\alpha)$$

$$= -r\zeta^{(r-1)}(s+1,\alpha) - s\zeta^{(r)}(s+1,\alpha)$$

so that we define

$$\frac{\partial}{\partial \alpha}\zeta^{(r)}(1,\alpha) = -r\zeta^{(r-1)}(2,\alpha) - \zeta^{(r)}(2,\alpha).$$

4) For $\alpha \neq 0,-1,-2,\ldots\ldots\ldots$, we have

$$\zeta(s,\alpha) = \frac{1}{s-1} + \sum_{r \geq 0} \gamma_r(\alpha)(s-1)^r,$$

where

$$\gamma_r(\alpha) = \gamma_r + \frac{1}{r!}(\zeta^{(r)}(s,\alpha) - \zeta^{(r)}(s))|_{s=1} \text{ and } \gamma_r = \gamma_r(1).$$

5) For $r \geq 0$, $\gamma_r(\alpha)$ is an analytic function of α except possibly for

$\alpha = 0,-1,-2,\ldots\ldots\ldots$, and we have

I)  $$\frac{\partial}{\partial \alpha}\zeta^{(r)}(0,\alpha) = -r!\gamma_{r-1}(\alpha), \text{ where } \gamma_{-1}(\alpha) = 1.$$

II) $$\frac{\partial}{\partial \alpha}\gamma_r(\alpha) = -\frac{\zeta^{(r-1)}(2,\alpha)}{(r-1)!} - \frac{\zeta^{(r)}(2,\alpha)}{r!} \text{ for } r \geq 1$$

and $$\frac{\partial}{\partial \alpha}\gamma_0(\alpha) = -\zeta(2,\alpha).$$

**Proof of Theorem:**

In author [2], it has been shown that if $\zeta_1(s,\alpha) = \zeta(s,\alpha) - \alpha^{-s}$, then for any complex $s$ and complex $\alpha$ with $|\alpha| < 1$,

$$\zeta_1(s,\alpha) = \zeta(s) + \sum_{n \geq 1} c_n(s)\frac{(-\alpha)^n}{n!},$$

where

$$c_n(s) = s(s+1)\ldots\ldots(s+n-1)\zeta(s+n).$$

Note that $B(s) = s\zeta(s+1)$ is an entire function of the complex variable $s$ with $B(0) = \lim_{s \to 0} s\zeta(s+1) = 1$. Thus for $n \geq 1$, $c_n(s)$ is an entire function of $s$.

More generally, for an integer $k \geq 1$ and for $0 < \alpha \leq 1$,

if $\zeta_k(s,\alpha) = \sum_{n \geq k}(n+\alpha)^{-s}$ for Re s>1; and its analytic continuation,

so that $$\zeta_k(s,\alpha) = \zeta(s,k+\alpha) = \zeta(s,\alpha) - \sum_{0 \leq n \leq k-1}(n+\alpha)^{-s},$$

then for any complex $s$ and for complex $\alpha$ with $|\alpha| < k$, we have

$$\zeta_k(s,\alpha) = \zeta_k(s) + \sum_{n \geq 1} c_{n,k}(s)\frac{(-\alpha)^n}{n!},$$

where $$c_{n,k}(s) = s(s+1)\ldots\ldots(s+n-1)\zeta_k(s+n)$$

and
$$\zeta_k(s) = \zeta(s) - \sum_{n \leq k-1} n^{-s},$$

with empty sum denoting zero.

Thus for $n \geq 1$, $c_{n,k}(s)$ is an entire function of $s$, as

$$s\zeta_k(s+1) = s(\zeta(s+1) - \sum_{n \leq k-1} n^{-s-1}) = s\zeta(s+1) - s \sum_{n \leq k-1} n^{-s-1}$$

is an entire function of $s$.

For any complex α with $|\alpha| < k$, and for any complex $s$, we have

$$\zeta(s,\alpha) - \zeta_k(s) = \sum_{0 \leq n \leq k-1} (n+\alpha)^{-s} + \sum_{n \geq 1} c_{n,k}(s) \frac{(-\alpha)^n}{n!}.$$

As the integer $k \geq 1$ is arbitrary, this defines $\zeta(s,\alpha)$ for any complex $s$ and α.

Note that if Re s≥0 and $s \neq 0$, then $\zeta(s,\alpha) - \zeta(s)$ is undefined for $\alpha = -n$ with $n = 0,1,2,\ldots,$.

If Re s<0, then $\zeta(s,\alpha)$ is defined for every $s$ and α.

Note that we have for $\alpha \neq 0,-1,-2,\ldots,$

$$(s-1)\zeta(s,\alpha) = (s-1)\ (\zeta(s) - \sum_{n \leq k-1} n^{-s}) + (s-1)\sum_{n \geq 1} c_{n,k}(s) \frac{(-\alpha)^n}{n!}.$$

Letting $s \to 1$, we have for $\alpha \neq 0,-1,\ldots,$ $\lim_{s \to 1}(s-1)\zeta(s,\alpha) = \lim_{s \to 1}(s-1)\zeta(s) = 1$.

Thus the function $B(s,\alpha) = (s-1)\zeta(s,\alpha)$ is well-defined for every $s$ and for $\alpha \neq 0,-1,-2,\ldots$ with $B(1,\alpha) = 1$ for $\alpha \neq 0,-1,-2,\ldots$.

As $B(1,\alpha)$ is an analytic function of α, we may define $B(1,\alpha) = 1$ for every complex number α.

We shall show that the series $\sum_{n\geq 1} c_{n,k}(s)\frac{(-\alpha)^n}{n!}$ as a series of analytic functions of α , is absolutely and uniformly convergent on every compact subset of the disc $|\alpha|<k$ for a fixed complex s . This will show that

$\zeta_k(s,\alpha) - \zeta_k(s) = \sum_{n\geq 1} c_{n,k}\frac{(-\alpha)^n}{n!}$ is an analytic function of α in the disc $|\alpha|<k$ .

Simultaneously , we shall show that the series $\sum_{n\geq 1} c_{n,k}(s)\frac{(-\alpha)^n}{n!}$ as a series of analytic functions of s , is absolutely and uniformly convergent on every compact subset of the complex s-plane , when $|\alpha|<k$ .

For $\sigma = \mathrm{Re}\ s > 1$ and for $|\alpha|<k$ , we have

$$\zeta_k(s,\alpha) = \sum_{n\geq k}(n+\alpha)^{-s} = \sum_{n\geq k} n^{-s}(1+\tfrac{\alpha}{n})^{-s} = \sum_{n\geq k} n^{-s}\{1 - s\tfrac{\alpha}{n} + \tfrac{s(s+1)}{2!}\tfrac{\alpha^2}{n^2} - \ldots\ldots\}$$

$$= \sum_{n\geq k} n^{-s} - \alpha s \sum_{n\geq k} n^{-s-1} + \tfrac{\alpha^2 s(s+1)}{2!}\sum_{n\geq k} n^{-s-2} - \ldots\ldots = \sum_{n\geq 0} c_{n,k}(s)\frac{(-\alpha)^n}{n!} .$$

We write $b_{n,k}(s) = \frac{(-1)^n}{n!} c_{n,k}(s)$ .

Thus for Re s>1 , we have $\zeta_k(s,\alpha) - \zeta_k(s) = \sum_{n\geq 1} b_{n,k}(s)\alpha^n$ .

Next we show the series $\sum_{n\geq 1} b_{n,k}(s)\alpha^n$ converges absolutely and uniformly in every compact subset of the half plane $\sigma>1$, when $|\alpha|<k$ ; and on every compact subset of the disc $|\alpha|<k$ , when $\sigma>1$. Hence consider $\beta$ such that $0<|\alpha|\leq \beta<k$ .

For $\sigma>1$, $|\sum_{n\geq 1} b_{n,k}(s)\alpha^n| \leq \sum_{n\geq 1}|b_{n,k}(s)||\alpha|^n \leq \sum_{n\geq 1}|b_{n,k}(s)|\beta^n$ .

Let $\beta' = \frac{\beta}{k}$ . Then $\beta'<1$.

Thus , we have $\qquad |\sum_{n\geq 1} b_{n,k}(s)\alpha^n| \leq \sum_{n\geq 1} k^n |b_{n,k}(s)|\beta'^n$ .

Next ,

$$k^n |b_{n,k}(s)| \leq \frac{k^n |s(s+1)....(s+n-1)| \cdot |\zeta_k(s+n)|}{n!} \leq \frac{|s|(|s|+1)....(|s|+n-1)k^n}{n!} \sum_{m \geq k} m^{-\sigma-n}$$

$$\leq |s|(|s|+1).....(|s|+n-1) \sum_{m \geq k} m^{-\sigma} \leq \frac{|s|(|s|+1)....(|s|+n-1)}{n!} \zeta(\sigma).$$

Thus 
$$|\sum_{n \geq 1} b_{n,k}(s)\alpha^n| \leq \zeta(\sigma) \sum_{n \geq 1} \frac{|s|(|s|+1)....(|s|+n-1)\beta'^n}{n!} \leq \zeta(\sigma)(1-\beta')^{-|s|}.$$

Thus in the disc $|\alpha| < k$, $\sum_{n \geq 1} b_{n,k}(s)\alpha^n$ defines an analytic function of α for a fixed s with Re s > 1, where $\zeta_k(s,\alpha) - \zeta_k(s) = \sum_{n \geq 1} b_{n,k}(s)\alpha^n$. This also proves that when Re s>1, $\sum_{n \geq 1} b_{n,k}(s)\alpha^n$ converges absolutely and uniformly as a series of analytic functions of s on every compact subset of the half plane Re s>1 for a fixed α with $|\alpha| < k$, as $|\sum_{n \geq 1} b_{n,k}(s)\alpha^n| \leq \zeta(\sigma)(1-\alpha')^{-|s|}$, where $\alpha' = \frac{|\alpha|}{k} < 1$. In the place of Re s>1, instead if we have Re s>-m for integral $m \geq 0$, the proof can be modified accordingly.

Thus $\zeta_k(s,\alpha) - \zeta_k(s) = \sum_{n \geq 1} c_{n,k}(s) \frac{(-\alpha)^n}{n!}$ defines an analytic function of s and we have for any integer $r \geq 1$,

$$\zeta_k^{(r)}(s,\alpha) - \zeta_k^{(r)}(s) = \sum_{n \geq 1} (c_{n,k}(s))^{(r)} \frac{(-\alpha)^n}{n!},$$

where the superscript $(r)$ denotes r-th order partial derivative with respect to s.

Next for any integer $r \geq 1$, and for any complex number s and for α with $|\alpha| < k$,

$$\zeta^{(r)}(s,\alpha) - \zeta_k^{(r)}(s) = \sum_{0 \leq n \leq k-1} (-1)^r (n+\alpha)^{-s} \log^r(n+\alpha)$$

$$+ \sum_{n \geq 1} \frac{(-\alpha)^n}{n!} \{s(s+1)..........(s+n-2) \cdot (s+n-1)\zeta_k(s+n)\}^{(r)},$$

where 
$$\zeta_k(s) = \zeta(s) - \sum_{n \leq k-1} n^{-s}$$

so that $(\zeta_k(s))^{(r)} = \zeta^{(r)}(s) - \sum_{n \leq k-1}(n^{-s})^{(r)} = \zeta^{(r)}(s) - \sum_{n \leq k-1}(-1)^r n^{-s} \log^r n$.

Thus for Re s<0 and for any α, $\zeta^{(r)}(s,\alpha)$ exists.

Thus for Re s<0 and for any α, $\zeta(s,\alpha)$ is an analytic function of the complex variable s.

If Re s≥0 then for $\alpha \neq 0,-1,-2,\ldots$, $\zeta^{(r)}(s,\alpha) - \zeta^{(r)}(s)$ exists and thus $\zeta(s,\alpha) - \zeta(s)$ is an analytic function of s for Re s≥0 and for $\alpha \neq 0,-1,-2,\ldots$

Thus if $s \neq 1$ and $\alpha \neq 0,-1,-2,\ldots$, $\zeta(s,\alpha)$ is an analytic function of s.

Next from the expression

$$\zeta(s,\alpha) = \zeta_k(s) + \sum_{0 \leq n \leq k-1}(n+\alpha)^{-s} + \sum_{n \geq 1} c_{n,k}(s) \frac{(-\alpha)^n}{n!} \text{ for } |\alpha| < k,$$

we have on differentiating with respect to α for |α|<k,

$$\frac{\partial \zeta(s,\alpha)}{\partial \alpha} = -s \sum_{0 \leq n \leq k-1}(n+\alpha)^{-s-1} - \sum_{n \geq 1} c_{n,k}(s) \frac{(-\alpha)^{n-1}}{(n-1)!}$$

$$= -s \sum_{0 \leq n \leq k-1}(n+\alpha)^{-s-1} - \sum_{n \geq 0} c_{n+1,k}(s) \frac{(-\alpha)^n}{n!}$$

$$= -s \sum_{0 \leq n \leq k-1}(n+\alpha)^{-s-1} - s \sum_{n \geq 0} c_{n,k}(s+1) \frac{(-\alpha)^n}{n!}$$

$$= -s \{ \sum_{0 \leq n \leq k-1}(n+\alpha)^{-s-1} + \sum_{n \geq 0} c_{n,k}(s+1) \frac{(-\alpha)^n}{n!} \} = -s\zeta(s+1,\alpha)$$

for $s \neq 1$ and for $\alpha \neq 0,-1,-2,\ldots$

Thus $\frac{\partial}{\partial \alpha} \zeta(s,\alpha) = -s\zeta(s+1,\alpha)$ for $\alpha \neq 0,-1,-2,\ldots$.

Note that $s\zeta(s+1,\alpha)$ exists for every s and for $\alpha \neq 0,-1,-2,\ldots$.

Thus $\zeta(s,\alpha)$ is an analytic function of the complex variable α for every s and for

$\alpha \neq 0, -1, -2, \ldots\ldots\ldots$ .

Thus if $s \neq 1$ and $\alpha \neq 0, -1, -2, \ldots\ldots$, $\zeta(s, \alpha)$ is an analytic function of $s$ and also of $\alpha$.

Next, we shall show that if $\alpha \neq 0, -1, -2, \ldots\ldots$,

then for $r \geq 1$, $\dfrac{\partial^r}{\partial s^r} \zeta(s, \alpha)$ is an analytic of $\alpha$ for every $s$ and

$$\frac{\partial}{\partial \alpha} \frac{\partial^r}{\partial s^r} \zeta(s, \alpha) = \frac{\partial^r}{\partial s^r} \frac{\partial}{\partial \alpha} \zeta(s, \alpha).$$

We have for $|\alpha| < k$ and for $r \geq 1$,

$$\zeta^{(r)}(s, \alpha) = \sum_{0 \leq n \leq k-1} (-1)^r (n+\alpha)^{-s} \log^r(n+\alpha) + \sum_{n \geq 0} \frac{(-\alpha)^n}{n!} \{s(s+1)\ldots\ldots(s+n-1)\zeta_k(s+n)\}^{(r)}.$$

Differentiating with respect to $\alpha$ for $|\alpha| < k$, we have

$$\frac{\partial}{\partial \alpha} \zeta^{(r)}(s, \alpha) = \sum_{0 \leq n \leq k-1} \frac{\partial}{\partial \alpha} \{(n+\alpha)^{-s}\}^{(r)} - \sum_{n \geq 1} \frac{(-\alpha)^{n-1}}{(n-1)!} \{s(s+1)\ldots\ldots(s+n-1)\zeta_k(s+n)\}^{(r)}.$$

$$= \sum_{0 \leq n \leq k-1} \frac{\partial^r}{\partial s^r} \frac{\partial}{\partial \alpha} (n+\alpha)^{-s} - \sum_{n \geq 0} \frac{(-\alpha)^n}{n!} \{s(s+1)\ldots\ldots(s+1+n-1) \cdot \zeta_k(s+1+n)\}^{(r)}$$

$$= \frac{\partial^r}{\partial s^r}(-s \sum_{0 \leq n \leq k-1} (n+\alpha)^{-s-1}) + \frac{\partial^r}{\partial s^r}\{-s \sum_{n \geq 0} \frac{(-\alpha)^n}{n!}(s+1)(s+2)\ldots\ldots(s+1+n-1) \cdot \zeta_k(s+1+n)\}$$

$$= \frac{\partial^r}{\partial s^r}(-s)\{\sum_{0 \leq n \leq k-1}(n+\alpha)^{-s-1} + \sum_{n \geq 0}\frac{(-\alpha)^n}{n!}(s+1)(s+2)\ldots\ldots(s+1+n-1)\zeta_k(s+1+n)\}$$

$$= \frac{\partial^r}{\partial s^r}(-s\zeta(s+1, \alpha)) = \frac{\partial^r}{\partial s^r} \frac{\partial}{\partial \alpha} \zeta(s, \alpha).$$

Thus $\dfrac{\partial^r}{\partial s^r} \zeta(s, \alpha)$ is an analytic function of $\alpha$ for $\alpha \neq 0, -1, -2, \ldots\ldots$ with

$$\frac{\partial}{\partial \alpha} \frac{\partial^r}{\partial s^r} \zeta(s, \alpha) = \frac{\partial^r}{\partial s^r} \frac{\partial}{\partial \alpha} \zeta(s, \alpha).$$

More generally , we have

$$\frac{\partial^{r_1}}{\partial \alpha^{r_1}} \frac{\partial^{r_2}}{\partial s^{r_2}} \zeta(s,\alpha) = \frac{\partial^{r_2}}{\partial s^{r_2}} \frac{\partial^{r_1}}{\partial \alpha^{r_1}} \zeta(s,\alpha)$$

for integral $r_1, r_2 \geq 0$.

Next we show that

$$\frac{\partial}{\partial \alpha} \frac{\partial^r}{\partial s^r} \zeta(s,\alpha) = -r\zeta^{(r-1)}(s+1,\alpha) - s\zeta^{(r)}(s+1,\alpha).$$

We have

$$\frac{\partial}{\partial \alpha} \frac{\partial^r}{\partial s^r} \zeta(s,\alpha) = \frac{\partial^r}{\partial s^r} \frac{\partial}{\partial \alpha} \zeta(s,\alpha) = \frac{\partial^r}{\partial s^r}(-s\zeta(s+1,\alpha))$$

$$= -\frac{\partial^r}{\partial s^r}(s\zeta(s+1,\alpha)) = -r\zeta^{(r-1)}(s+1,\alpha) - s\zeta^{(r)}(s+1,\alpha).$$

Next, we know that for any complex s and for α with $0 < \alpha \leq 1$, we have

$$\zeta(s,\alpha) = \frac{1}{s-1} + \sum_{r \geq 0} \gamma_r(\alpha)(s-1)^r,$$

where

$$\gamma_r(\alpha) = \frac{1}{r!}\{\zeta(s,\alpha) - \frac{1}{s-1}\}^{(r)}_{s=1}$$

$$= \frac{1}{r!}\{(\zeta(s,\alpha) - \zeta(s)) + (\zeta(s) - \frac{1}{s-1})\}^{(r)}_{s=1} = \frac{1}{r!}(\zeta(s,\alpha) - \zeta(s))^{(r)}_{s=1} + \gamma_r,$$

where

$$\zeta(s) = \frac{1}{s-1} + \sum_{r \geq 0} \gamma_r(s-1)^r.$$

Thus

$$\gamma_r(\alpha) - \gamma_r = \frac{1}{r!}(\zeta(s,\alpha) - \zeta(s))^{(r)}_{s=1}.$$

But we know that

$$\zeta^{(r)}(s,\alpha) - \zeta^{(r)}(s)|_{s=1}$$

is an analytic function of α.

Thus $\gamma_r(\alpha)$ is an analytic function of α for $\alpha \neq 0, -1, -2, \ldots\ldots\ldots\ldots$.

Next we evaluate $\dfrac{\partial}{\partial \alpha}\gamma_r(\alpha)$.

We have $\dfrac{\partial}{\partial \alpha}\gamma_r(\alpha)$

$$= \frac{1}{r!}\frac{\partial}{\partial \alpha}(\zeta(s,\alpha)-\zeta(s))^{(r)}_{s=1} = \frac{1}{r!}\frac{\partial}{\partial \alpha}\frac{\partial^r}{\partial s^r}(\zeta(s,\alpha)-\zeta(s))|_{s=1} = \frac{1}{r!}\frac{\partial^r}{\partial s^r}\frac{\partial}{\partial \alpha}(\zeta(s,\alpha)-\zeta(s))|_{s=1}$$

$$= \frac{1}{r!}\frac{\partial^r}{\partial s^r}\frac{\partial}{\partial \alpha}\zeta(s,\alpha)|_{s=1} = \frac{1}{r!}\frac{\partial^r}{\partial s^r}(-s\zeta(s+1,\alpha))|_{s=1} = \frac{1}{r!}\{-r\zeta^{(r-1)}(s+1,\alpha) - s\zeta^{(r)}(s+1,\alpha)\}_{s=1}$$

$$= -\zeta(2,\alpha) \text{ for } r=0$$

$$= -\frac{\zeta^{(r-1)}(2,\alpha)}{(r-1)!} - \frac{\zeta^{(r)}(2,\alpha)}{r!} \text{ , if } r \geq 1.$$

and

Next we evaluate

$$\frac{\partial}{\partial \alpha}\frac{\partial^r}{\partial s^r}\zeta(s,\alpha) \text{ at } s=0.$$

We have

$$\frac{\partial}{\partial \alpha}\frac{\partial^r}{\partial s^r}\zeta(s,\alpha)|_{s=0} = \frac{\partial^r}{\partial s^r}\frac{\partial}{\partial \alpha}\zeta(s,\alpha)|_{s=0} = \frac{\partial^r}{\partial s^r}(-s\zeta(s+1,\alpha))|_{s=0}.$$

But we know that

$$s\zeta(s+1,\alpha) = 1 + \sum_{r\geq 1}\gamma_{r-1}(\alpha)s^r$$

so that $\gamma_{r-1}(\alpha) = \dfrac{1}{r!}(s\zeta(s+1,\alpha))^{(r)}_{s=0}$.

Thus
$$\frac{\partial}{\partial \alpha} \frac{\partial^r}{\partial s^r} \zeta(s,\alpha)|_{s=0} = -r!\gamma_{r-1}(\alpha) \ .$$